\documentclass{amsart}
\usepackage{amscd}
\usepackage{amsmath}
\usepackage{amssymb}
\usepackage{amsthm}
\newtheorem{theorem}{Theorem}[section]
\newtheorem{lemma}[theorem]{Lemma}
\newtheorem{proposition}[theorem]{Proposition}
\newtheorem{corollary}[theorem]{Corollary}

\newtheorem{problem}[theorem]{Problem}
\theoremstyle{definition}

\newtheorem{example}[theorem]{Example}

\newtheorem*{acknowledgements}{Acknowledgments}
\theoremstyle{remark}

\numberwithin{equation}{section}
\DeclareMathOperator{\Tor}{Tor}
\DeclareMathOperator{\pd}{pd}

\DeclareMathOperator{\height}{height}

\DeclareMathOperator{\indeg}{indeg}

\DeclareMathOperator{\reg}{reg}

\DeclareMathOperator{\ara}{ara}
\DeclareMathOperator{\chara}{char}

\DeclareMathOperator{\cone}{co}

\begin{document}
%%%%%%%%%%%%%%%%%%%%%%%%%%%%%%%%%%%%%%%%
%% title
%%%%%%%%%%%%%%%%%%%%%%%%%%%%%%%%%%%%%%%%
\title[Arithmetical rank of Cohen--Macaulay monomial ideals of height two]
{Arithmetical rank of Cohen--Macaulay squarefree monomial ideals of height two}
%%%%%%%%%%%%%%%%%%%%%%%%%%%%%%%%%%%%%%%%
%%%%%%%%%%%%%%%%%%%%%%%%%%%%%%%%%%%%%%%%
%% Information for authors
%%%%%%%%%%%%%%%%%%%%%%%%%%%%%%%%%%%%%%%%
\author[Kyouko Kimura]{Kyouko Kimura}
%\author[K. Kimura]{Kyouko Kimura}
\address{Department of Pure and Applied Mathematics, 
Graduate School of Information Science and Technology, Osaka University, 
Toyonaka, Osaka 560-0043, Japan}
\email{kimura@math.sci.osaka-u.ac.jp}

%%%%%%%%%%%%%%%%%%%%%%%%%%%%%%%%%%%%%%%%
%% General info
%%%%%%%%%%%%%%%%%%%%%%%%%%%%%%%%%%%%%%%%
\date{\today}
%\dedicatory{}
%%%%%%%%%%%%%%%%%%%%%%%%%%%%%%%%%%%%%%%%
%%%%%%%%%%%%%%%%%%%%%%%%%%%%%%%%%%%%%%%%
%% Abstract
%%%%%%%%%%%%%%%%%%%%%%%%%%%%%%%%%%%%%%%%
\begin{abstract}
In this paper, we prove that a squarefree monomial ideal of height $2$ 
whose quotient ring is Cohen--Macaulay is set-theoretic complete intersection. 
\end{abstract}

%%%%%%%%%%%%%%%%%%%%%%%%%%%%%%%%%%%%%%%%
\maketitle
%\thispagestyle{empty}
%%%%%%%%%%%%%%%%%%%%%%%%%%%%%%%%%%%%%%%%
\section{Introduction}
Let $R$ be a polynomial ring over a field $K$. 
Let $I$ be a squarefree monomial ideal of $R$ 
and $G(I)$ the minimal set of monomial generators of $I$. 
The \textit{arithmetical rank} of $I$ is defined by the minimum number $r$ 
of elements $a_1, \ldots, a_r \in R$ such that 
\begin{equation}
  \label{eq:up_to_radical}
  \sqrt{(a_1, \ldots, a_r)} = \sqrt{I}. 
\end{equation}
We denote it by $\ara I$. 
When (\ref{eq:up_to_radical}) holds, we say that 
\textit{$a_1, \ldots, a_r$ generate $I$ up to radical}. 
By Krull's principal ideal theorem, we have $\height I \leq \ara I$. 
When the equality holds, 
we say that $I$ is \textit{set-theoretic complete intersection}. 
Moreover, Lyubeznik \cite{Ly83} proved that 
for a squarefree monomial ideal $I$, 
the projective dimension of $R/I$ over $R$, 
denoted by $\pd_R R/I$ (or $\pd R/I$ if there is no confusion), 
provides a better lower bound of the arithmetical rank of $I$. 
Many authors involving Barile 
\cite{Barile96, Barile08-1, Barile08-2, Barile0606}, 
Barile and Terai \cite{BariTera08, BariTera09}, 
Ene, Olteanu and Terai \cite{EOT}, Kummini \cite{Kummini}, 
Schmitt and Vogel \cite{SchmVo}, 
Terai and Yoshida with the author \cite{KTYdev1, KTYdev2}, 
investigated when $\ara I = \pd_R R/I$ holds. 

\par
In this paper, we prove the following theorem: 
\begin{theorem}[See Theorem \ref{claim:h2CM}]
  \label{claim:h2CM_intro}
  Let $I$ be a squarefree monomial ideal of $R$ of height $2$. 
  Suppose that $R/I$ is Cohen--Macaulay. Then 
  \begin{displaymath}
    \ara I = \pd_R R/I = 2. 
  \end{displaymath}
  In particular, $I$ is set-theoretic complete intersection. 
\end{theorem}
That is, ideals as in Theorem \ref{claim:h2CM_intro} 
are generated by $2$ elements up to radical. 
Note that the equality $\ara I = \pd_R R/I$ does not always hold 
for Cohen--Macaulay squarefree monomial ideals $I$ of height $3$ 
(when $\chara K \neq 2$) as founded by Yan \cite{Yan}, 
Terai and Yoshida with the author \cite{KTYdev2}. 

\par
We explain the organization of this paper. 
First in Section $2$, we state the motivated problem 
of this paper (Problem \ref{prob:motivation}), 
which corresponds to Alexander dual of the results 
in Barile and Terai \cite{BariTera08}.  
Partial answers for this problem are given in Section 3 
(Propositions \ref{claim:ara} and \ref{claim:ara=2}). 
In particular, Proposition \ref{claim:ara=2} plays the key role on the 
proof of Theorem \ref{claim:h2CM_intro}, which is given in Section 4. 
%In Section 4, we prove Theorem \ref{h2CM_intro}. 

\par
The main result of Barile and Terai \cite[Theorem 1]{BariTera08}, 
which is our motivated paper, required the assumption that 
$K$ is algebraically closed. 
At the last of this paper, in Section $5$, 
we give an improvement proof of that result. 
Consequently, we can remove the assumption on $K$. 

\section{Preliminaries and the motivated problem}
In this section, we state the motivated problem of this paper. 
As before, we recall some definitions and properties of 
simplicial complexes and Stanley--Reisner ideals, 
especially, Alexander duality.  
For more detail, we refer to \cite[Section 5]{Bruns-Herzog}, \cite{Terai07}. 

\par
Let $I$ be a squarefree monomial ideal of a polynomial ring $R$ 
over a field $K$. 
The \textit{graded Betti number} of $R/I$ is defined by 
${\beta}_{i, j} (R/I) = \dim_K [ \Tor_i^R (R/I, K)]_j$. 
The \textit{initial degree}, 
the \textit{$($Castelnuovo--Mumford$)$ regularity} 
of $I$ are defined by 
\begin{displaymath}
  \indeg I = \min \{ j \; : \; {\beta}_{1, j} (R/I) \neq 0 \}, 
  \quad \reg I = \max \{ j-i+1 \; : \; {\beta}_{i, j} (R/I) \neq 0 \}, 
\end{displaymath}
respectively. 
In general, the inequality $\reg I \geq \indeg I$ holds. 
When $\reg I = \indeg I = k$, 
we say that $I$ has a \textit{$k$-linear resolution}. 

\par
Let $X = \{ x_1, x_2, \ldots, x_n \}$ be a set of indeterminates 
over a field $K$. 
A \textit{simplicial complex} $\Delta$ on the vertex set $X$ is 
a collection of subsets of $X$ with the properties 
(i) $\{ x_i \} \in \Delta$ for all $x_i \in X$; 
(ii) $F \in \Delta$ and $G \subset F$ imply $G \in \Delta$. 
If $\Delta$ consists of all subsets of $X$, then $\Delta$ is called a 
\textit{simplex}. 
An element of $\Delta$ is called a \textit{face} of $\Delta$. 
A maximal face of $\Delta$ with respect to inclusion is called a 
\textit{facet} of $\Delta$. 
The \textit{dimension} of $\Delta$ is defined by 
$\dim \Delta = \max \{ |F| - 1 : F \in \Delta \}$, 
where $|F|$ denotes the cardinality of $F$. 
The \textit{Alexander dual complex} ${\Delta}^{\ast}$ is defined by 
${\Delta}^{\ast} = \{ F \subset X : X \setminus F \notin \Delta \}$, 
which is also a simplicial complex. 
If $\dim \Delta < n - 2$, 
then the vertex set of ${\Delta}^{\ast}$ coincides with $X$. 
When this is the case, ${\Delta}^{\ast \ast} = \Delta$. 

\par
For a simplicial complex $\Delta$ on the vertex set 
$X = \{ x_1, x_2, \ldots, x_n \}$, we associate a squarefree monomial 
ideal $I_{\Delta}$ of $K[X] = K[x_1, x_2, \ldots, x_n]$ as follows: 
\begin{displaymath}
  I_{\Delta} = \big( x_{i_1} \cdots x_{i_s} \; : \; 
                     1 \leq i_1 < \cdots < i_s \leq n, \  
                     \{ x_{i_1}, \ldots, x_{i_s} \} \notin \Delta \big), 
\end{displaymath}
which is called the \textit{Stanley--Reisner ideal} of $\Delta$. 
The quotient ring $K[{\Delta}] = K[X]/{I_{\Delta}}$ is called the 
\textit{Stanley--Reisner ring} of $\Delta$. 
The minimal prime decomposition of $I_{\Delta}$ is given by 
\begin{equation}
  \label{eq:prim_decomp}
  I_{\Delta} = \bigcap_{F \in \Delta \; : \; \text{facet}} 
                   P_{F}, 
\end{equation}
where $P_{F} = (x_i : x_i \in X \setminus F)$. 

\par
On the other hand, it is well-known that for a squarefree monomial ideal 
$I$ of $R = K[X]$ with $\indeg I \geq 2$, 
there exists a simplicial complex $\Delta$ on $X$ such that $I = I_{\Delta}$. 
Assume that $\height I \geq 2$. 
%Then since the Alexander dual complex ${\Delta}^{\ast}$ 
%is also a simplicial complex on the same vertex set $X$, 
Then since $\dim \Delta < n-2$, 
we can consider the ideal $I^{\ast} = I_{{\Delta}^{\ast}}$ of $R$, 
which is called the \textit{Alexander dual ideal} of $I = I_{\Delta}$.  
Since ${\Delta}^{\ast \ast} = \Delta$, we have $I^{\ast \ast} = I$. 
The minimal set of monomial generators of 
$I^{\ast} = I_{{\Delta}^{\ast}}$ is given by 
\begin{equation}
  \label{eq:Adual_gen}
  G (I^{\ast}) = G(I_{{\Delta}^{\ast}}) 
  = \{ m_{X \setminus F} \; : \; 
        \text{$F \in \Delta$ is a facet of $\Delta$} \},  
\end{equation}
where $m_{X \setminus F} = \prod_{x_i \in X \setminus F} x_i$. 
Then it is easy to see that $\indeg I^{\ast} = \height I$ 
by (\ref{eq:prim_decomp}), (\ref{eq:Adual_gen}). 
Moreover, 
%Terai proved that $\reg I^{\ast} = \pd_R R/I$ 
%(see Terai \cite[Corollary 1.6]{Terai07}) and 
Eagon and Reiner \cite[Theorem 3]{Ea-Re} proved that 
$I$ has a linear resolution if and only if $R/{I^{\ast}}$ is Cohen--Macaulay. 
%(see \cite[Theorem 1.3]{Terai07}). 

\par
\bigskip

\par
Now we state our motivated problem. 

\par
Let $\Delta$ be a simplicial complex 
on the vertex set $X = \{ x_1, x_2, \ldots, x_n \}$. 
Let $x_0$ be a new indeterminate and $F$ a face of $\Delta$. 
A \textit{cone from $x_0$ over $F$}, 
denoted by $\cone_{x_0} F$, is the simplex on the vertex set 
$F \cup \{ x_0 \}$. 
%We denote it by $\cone_{x_0} F$. 
Then ${\Delta}' := \Delta \cup \cone_{x_0} F$ is a simplicial complex 
on the vertex set $X' := X \cup \{ x_0 \}$. 
Barile and Terai \cite{BariTera08} investigated some relations between 
arithmetical ranks of $I_{\Delta}$ and $I_{{\Delta}'}$ 
(\cite[Theorem 1]{BariTera08}). 
Moreover, they proved that 
if $\ara I_{\Delta} = \pd K[{\Delta}]$ holds, 
then $\ara I_{{\Delta}'} = \pd K[{\Delta}']$ also holds 
(\cite[Theorem 2]{BariTera08}). 
As its corollary, they proved that if a squarefree monomial ideal 
$I \subset R$ has a $2$-linear resolution, 
then $\ara I = \pd_R R/I$ holds. 
(This result was first proved by Morales \cite{Mo} on the different way.) 

\par
We consider the following problem which corresponding to 
Alexander dual of their results: 
\begin{problem}
  \label{prob:motivation}
  Let $\Delta$ be a simplicial complex on the vertex set 
  $X = \{ x_1, x_2, \ldots, x_n \}$ 
  with $\dim \Delta < n - 2$. 
  Let $F$ be an arbitrary face of ${\Delta}^{\ast}$ and $x_0$ a new vertex. 
  Set $X' = X \cup \{ x_0 \}$, ${\Gamma} = {\Delta}^{\ast}$, 
  ${\Gamma}' = {\Gamma} \cup \cone_{x_0} F$, 
  and ${\Delta}' = ({\Gamma}')^{\ast}$. 

  \par
  Are there any relations between arithmetical ranks 
  of $I_{\Delta}$ and $I_{{\Delta}'}$? 
  In particular, if $\ara I_{\Delta} = \pd K[{\Delta}]$ holds, 
  then does $\ara I_{{\Delta}'} = \pd K[{\Delta}']$ hold? 
\end{problem}

\par
Here, we compare the projective dimension of $K[{\Delta}']$ 
with that of $K[{\Delta}]$. 

\begin{lemma}
  \label{claim:proj_dim}
  Let $\Delta$ and ${\Delta}'$ be simplicial complexes 
  as in Problem \ref{prob:motivation}. 
  Then 
  \begin{displaymath}
    \pd K[{\Delta}'] = \pd K[{\Delta}]. 
  \end{displaymath}
\end{lemma}
\begin{proof}
  Set $R = K[X]$, $R' = K[X']$, $I = I_{\Delta}$, $I' = I_{{\Delta}'}$, 
  and $G(I) = \{ m_1, \ldots, m_{\mu} \}$. 
  Then 
  \begin{displaymath}
      I_{\Gamma} = I^{\ast} 
      = P_{G_1} \cap \cdots \cap P_{G_{\mu}} \subset R, 
  \end{displaymath}
  where 
%  $P_i := P_{G_i}$ is a prime ideal generated 
%  by all variables $x_j \in X$ which divides $m_i$, 
  $G_1, \ldots, G_{\mu}$ are all facets of 
  $\Gamma = {\Delta}^{\ast}$ and $m_j = \prod_{x_i \in P_{G_j}} x_i$. 
  We may assume $F \subset G_1$ without loss of the generality. Then 
  \begin{displaymath}
    I_{{\Gamma}'} = P_{F \cup \{ x_0 \} } \cap (P_{G_1} R' + (x_0)) 
                    \cap \cdots \cap (P_{G_{\mu}} R' + (x_0)) \subset R'. 
  \end{displaymath}
  Hence 
  \begin{displaymath}
    I' = (m_0, x_0 m_1, \ldots, x_0 m_{\mu}) R' = m_0 R' + x_0 I R', 
  \end{displaymath}
  where $m_0 = \prod_{x_i \in P_{F \cup \{ x_0 \} }} x_i$. 
  Note that $m_0$ is divisible by $m_1$ 
  since $P_{G_1} R' \subset P_{F \cup \{ x_0 \}}$. 
  
  \par
  Let us consider the short exact sequence 
  \begin{equation}
    \label{eq:exact}
    0 \longrightarrow R'/{m_0 R' \cap x_0 I R'} 
      \longrightarrow R'/{m_0 R'} \oplus R'/{x_0 I R'} 
      \longrightarrow R'/I' \longrightarrow 0. 
  \end{equation}
  Note that $m_0 R' \cap x_0 I R' = x_0 m_0 R'$. 
  Since $\pd_{R'} R'/I = \pd_R R/I \geq \height I \geq 2$, 
  $\pd_{R'} R'/{x_0 m_0 R'} = \pd_{R'} R'/{m_0 R'} = 1$, 
  the long exact sequence obtained by applying $\Tor^{R'} (-, K)$ 
  to (\ref{eq:exact}) yields 
  \begin{displaymath}
    \pd K[{\Delta}'] = \pd_{R'} R'/{I'} = \pd_{R} R/{I} = \pd K[{\Delta}], 
  \end{displaymath}
  as desired. 
\end{proof}

\section{Partial answers for Problem \ref{prob:motivation}}
In this section, we give partial answers for Problem \ref{prob:motivation}. 
Throughout of this section, 
we use the notations as in Problem \ref{prob:motivation}. 

\par
First, we show a relation between arithmetical ranks 
of $I_{\Delta}$ and $I_{{\Delta}'}$. 
\begin{proposition}
  \label{claim:ara}
  Let $\Delta$ and ${\Delta}'$ be simplicial complexes 
  as in Problem \ref{prob:motivation}. 
  Then 
  \begin{displaymath}
    \ara I_{{\Delta}'} \leq \ara I_{\Delta} + 1. 
  \end{displaymath}

  \par
  In particular, if $\ara I_{\Delta} = \pd K[{\Delta}]$ holds, then 
  $\ara I_{{\Delta}'}$ coincides with either $\pd K[{\Delta}']$ 
  or $\pd K[{\Delta}'] + 1$. 
\end{proposition}
\begin{proof}
  Set $R = K[X]$, $R' = K[X']$, $I = I_{\Delta}$, $I' = I_{{\Delta}'}$, 
  and $G(I) = \{ m_1, \ldots, m_{\mu} \}$. Then 
  \begin{displaymath}
    I' = (m_0, x_0 m_1, \ldots, x_0 m_{\mu}) R', 
  \end{displaymath}
  where $m_0 = \prod_{x_i \in P_{F \cup \{ x_0 \} }} x_i$. 

  \par
  Put $h = \ara I$ and let $q_1, \ldots, q_h$ be elements of $R$ 
  which generate $I$ up to radical. 
  Then $x_0 q_1, \ldots, x_0 q_h$ generate $(x_0 m_1, \ldots, x_0 m_{\mu})$ 
  up to radical. 
  This implies that $m_0, x_0 q_1, \ldots, x_0 q_h$ generate 
  $I'$ up to  radical. 
  Therefore we have $\ara I' \leq h+1$. 

  \par
  Then the latter claim immediately follows from Lemma \ref{claim:proj_dim} 
  and the inequality $\ara I_{{\Delta}'} \geq \pd K[{\Delta}']$. 
\end{proof}

\par
Next, we give a partial answer for the second question 
of Problem \ref{prob:motivation}. 
\begin{proposition}
  \label{claim:ara=2}
  Let $\Delta$ and ${\Delta}'$ be simplicial complexes 
  as in Problem \ref{prob:motivation}. 
  Suppose that $\ara I_{\Delta} = \pd K[{\Delta}] = 2$. Then 
  \begin{displaymath}
    \ara I_{{\Delta}'} = \pd K[{\Delta}'] = 2. 
  \end{displaymath}
\end{proposition}

\par
In the study of the arithmetical rank, the technique based on 
linear algebraic consideration has been developed 
by Barile \cite{Barile08-1}, 
Barile and Terai \cite{BariTera08} (see also \cite{BariTera09}). 
Our proof of this proposition also goes along this current. 

\begin{proof}
  By Lemma \ref{claim:proj_dim}, we have 
  $\pd K[{\Delta}'] = \pd K[{\Delta}] = 2$. 
  Therefore it suffices to prove that $\ara I_{{\Delta}'} \leq 2$. 

  \par
  Set $R = K[X]$, $R' = K[X']$, $I = I_{\Delta}$, $I' = I_{{\Delta}'}$, 
  and $G(I) = \{ m_1, \ldots, m_{\mu} \}$. 
  Then 
  \begin{displaymath}
    I' = (m_0, x_0 m_1, \ldots, x_0 m_{\mu}) R', 
  \end{displaymath}
  where $m_0 = \prod_{x_i \in P_{F \cup \{ x_0 \} }} x_i$. 
  Let $G$ be a facet of $\Gamma = {\Delta}^{\ast}$ containing $F$. 
  We may assume $m_1 = \prod_{x_i \in P_G} x_i$. 
  Then $m_0$ is divisible by $m_1$ 
  since $P_G R' \subset P_{F \cup \{ x_0 \} }$. 

  \par
  Let $q_1, q_2$ be elements of $R$ which generate $I$ up to radical. 
  Note that $q_1, q_2 \in I$ because $I$ is a squarefree monomial ideal. 
  By $m_i \in \sqrt{(q_1, q_2)}$, there exists some integer ${\ell}_i \geq 1$ 
  such that $m_i^{{\ell}_i} \in (q_1, q_2)$. 
  Then we can write 
  \begin{displaymath}
    m_i^{{\ell}_i} = a_{i 1} q_1 + a_{i 2} q_2, \quad i = 1, \dots, \mu, 
  \end{displaymath}
  where $a_{i 1}, a_{i 2} \in R$. 
  Set $A = (a_{i j})_{i=1, \ldots, \mu; j=1, 2}$. Then 
  \begin{displaymath}
    \left(
    \begin{matrix}
      m_1^{{\ell}_1} \\
      \vdots \\
      m_{\mu}^{{\ell}_{\mu}} 
    \end{matrix}
    \right) = A \left( 
    \begin{matrix}
      q_1 \\
      q_2 
    \end{matrix}
    \right). 
  \end{displaymath}

  \par
  Set 
  \begin{displaymath}
    J' = (x_0 q_1 - a_{1 2} m_0, x_0 q_2 + a_{1 1} m_0) R'. 
  \end{displaymath}
  We prove $\sqrt{J'} = I'$. 
  Since $x_0 q_1 - a_{1 2} m_0, x_0 q_2 + a_{1 1} m_0 \in I'$, 
  we have $\sqrt{J'} \subset I'$. 
  We prove the opposite inclusion. 

  \par
  Since 
  \begin{displaymath}
    A \left( 
    \begin{matrix}
      x_0 q_1 - a_{1 2} m_0 \\
      x_0 q_2 + a_{1 1} m_0
    \end{matrix}
    \right) = \left( 
    \begin{matrix}
      x_0 m_1^{{\ell}_1} + f_1 m_0 \\
      \vdots \\
      x_0 m_{\mu}^{{\ell}_{\mu}} + f_{\mu} m_0 
    \end{matrix}
    \right), 
    \quad \text{where} \quad 
    f_i = a_{1 1} a_{i 2} - a_{1 2} a_{i 1}, 
  \end{displaymath}
  we have $x_0 m_i^{{\ell}_i} + f_i m_0 \in J'$ for $i=1, \ldots, \mu$. 
  Note that $f_1 = a_{1 1} a_{1 2} - a_{1 2} a_{1 1} = 0$. 
  Thus $x_0 m_1^{{\ell}_1} \in J'$, that is, $x_0 m_1 \in \sqrt{J'}$. 
  Since $m_1$ divides $m_0$, multiplying $x_0 m_i^{{\ell}_i} + f_i m_0 \in J'$ 
  by $x_0$ implies $x_0^2 m_i^{{\ell}_i} \in \sqrt{J'}$, 
  that is, $x_0 m_i \in \sqrt{J'}$. 

  \par
  Here, recall that $q_1, q_2 \in I = (m_1, \ldots, m_{\mu})$. 
%  Thus we can write 
%  $q_j = \sum_{i=1}^{\mu} b_{j i} m_i$ for $i=1, 2$. 
  Thus $x_0 q_1, x_0 q_2 \in \sqrt{J'}$. 
  Consequently, we have $a_{1 1} m_0, a_{1 2} m_0 \in \sqrt{J'}$. 
  By $m_1^{{\ell}_1} = a_{1 1} q_1 + a_{1 2} q_2$, we have 
  \begin{displaymath}
    m_0 m_1^{{\ell}_1} = m_0 (a_{1 1} q_1 + a_{1 2} q_2) 
      =(a_{1 1} m_0) q_1 + (a_{1 2} m_0) q_2 \in \sqrt{J'}. 
  \end{displaymath}
  This implies $m_0 \in \sqrt{J'}$ since $m_0$ is divisible by $m_1$. 
  Therefore $\sqrt{J'} \supset I'$ holds, as required. 
\end{proof}

\begin{example}
  \label{exam:AdualLine}
  Let $\Delta$ be the simplicial complex 
  on the vertex set $\{ x_1, x_2, x_3, x_4 \}$ whose facets are 
  $\{ x_1, x_3 \}, \{ x_2, x_3 \}, \{ x_2, x_4 \}$. 
  Then 
  \begin{displaymath}
%    \begin{aligned}
      I = I_{\Delta} 
        = (x_2, x_4) \cap (x_1, x_4) \cap (x_1, x_3) 
        = (x_1 x_2, x_1 x_4, x_3 x_4). 
%    \end{aligned}
  \end{displaymath}
  The Alexander dual complex $\Gamma$ of $\Delta$ has facets 
  $\{ x_3, x_4 \}, \{ x_2, x_3 \}, \{ x_1, x_2 \}$, that is, 
  $\Gamma$ is a line segment with $4$ vertices. 
  Take the face $F = \{ x_4 \} \in \Gamma$ and a new vertex $x_5 := x_0$. 
  Then ${\Gamma}' = \Gamma \cup \cone_{x_5} \{ x_4 \}$ 
  is a line segment with $5$ vertices and 
  \begin{displaymath}
    I_{{\Gamma}'} = (x_1, x_2, x_3) \cap (x_1, x_2, x_5) \cap (x_1, x_4, x_5) 
      \cap (x_3, x_4, x_5). 
  \end{displaymath}
  Thus $I' = I_{({\Gamma}')^{\ast}}$ is generated by 
  \begin{displaymath}
    x_1 x_2 x_3,\   x_1 x_2 x_5,\   x_1 x_4 x_5,\   x_3 x_4 x_5. 
  \end{displaymath}

  \par
  In this case, $m_0 = x_1 x_2 x_3$ and $m_1 = x_1 x_2$. 
%  We set $m_2 = x_1 x_4$ and $m_3 = x_3 x_4$. 
  By the result of Schmitt and Vogel \cite[Lemma, p.\  249]{SchmVo}, 
  it is easy to see that 
  the following two elements $q_1, q_2$ generate $I$ up to radical: 
  \begin{displaymath}
    q_1 = x_1 x_4, \qquad q_2 = x_1 x_2 + x_3 x_4. 
  \end{displaymath}
  Then
  \begin{equation}
    \label{eq:AdualLine}
    m_1^2 = - x_2 x_3 q_1 + x_1 x_2 q_2. 
  \end{equation}
  By Proposition \ref{claim:ara=2}, the following two elements 
  $q_1', q_2'$ generate $I'$ up to radical: 
  \begin{displaymath}
    \begin{aligned}
      q_1' &= x_5 q_1 - x_1 x_2 m_0 
            = x_1 x_4 x_5 -x_1^2 x_2^2 x_3, \\
      q_2' &= x_5 q_2 - x_2 x_3 m_0 
            = x_1 x_2 x_5 + x_3 x_4 x_5 - x_1 x_2^2 x_3^2. 
    \end{aligned}
  \end{displaymath}
\end{example}

\section{Proof of the main theorem}
In this section, we prove the following theorem, 
which is the main result in this paper. 
\begin{theorem}
  \label{claim:h2CM}
  Let $I$ be a squarefree monomial ideal of $R = K[X]$ of height $2$. 
  Suppose that $R/I$ is Cohen--Macaulay. 
  Then 
  \begin{displaymath}
    \ara I = \pd_R R/I = \height I = 2. 
  \end{displaymath}
  In particular, $I$ is set-theoretic complete intersection. 
\end{theorem}

\par
The Alexander dual of ideals satisfying the assumptions of this theorem 
have a $2$-linear resolution. 
To study these ideals, we recall the definition of the generalized tree. 

\par
We say that a simplicial complex is a \textit{generalized tree} 
if it can be obtained by the following recursive procedure: 
(i) a simplex is a generalized tree; (ii) if $\Delta$ is a generalized tree, 
then $\Delta \cup \cone_{x_0} F$ is also a generalized tree 
for any $F \in \Delta$ and for any new vertex $x_0$. 
Then a Stanley--Reisner ideal $I_{\Delta}$ which has a $2$-linear resolution 
is characterized as the following lemma. 
\begin{lemma}[{See Barile and Terai \cite[Lemma 2]{BariTera08}}]
  \label{claim:2-linear_chara}
  Let $\Delta$ be a simplicial complex which is not a simplex. 
  Then $I_{\Delta}$ has a $2$-linear resolution 
  if and only if $\Delta$ is a generalized tree. 
\end{lemma}

\par
Now we prove Theorem \ref{claim:h2CM}. 
The proof is done as an application of Proposition \ref{claim:ara=2}. 
\begin{proof}[Proof of Theorem \ref{claim:h2CM}]
  Since $R/I$ is Cohen--Macaulay, we have $\pd_R R/I = \height I = 2$. 
  First, we note that when $\mu (I) \leq \pd_R R/I + 1$, it is known that 
  $\ara I = \pd_R R/I$ holds; see e.g., \cite[Theorem 2.1]{KTYdev1}. 
  Thus in our situation, 
  $\ara I = \pd_R R/I = 2$ holds if $\mu (I) \leq 3$. 

  \par
  If $\indeg I = 1$, then $I$ is of the form $(x_1, m_2)$ 
  by the assumptions of $I$. In this case, $\ara I = 2$ trivially holds. 

  \par
  Assume that $\indeg I \geq 2$. 
  We proceed the proof by induction 
  on the number $|X|$ of variables. 
  The minimum number $|X|$ in which there exists an ideal $I$ 
  satisfying our assumption is $3$ and such an ideal is of the form 
  \begin{displaymath}
    I = (x_1, x_2) \cap (x_1, x_3) \cap (x_2, x_3) 
      = (x_1 x_2, x_1 x_3, x_2 x_3). 
  \end{displaymath}
  Then since $\mu (I) = 3$, we have $\ara I = \pd_R R/I = 2$. 

  \par
  Now assume $|X| > 3$. Since $I^{\ast} = I_{\Gamma}$ has 
  a $2$-linear resolution, $\Gamma$ is a generalized tree 
  by Lemma \ref{claim:2-linear_chara}, and there exist a vertex $x \in X$, 
  a generalized tree $\overline{\Gamma}$ on the vertex set 
  $X \setminus \{ x \}$ and a face $F \in \overline{\Gamma}$ such that 
  $\Gamma = \overline{\Gamma} \cup \cone_x F$ by the definition of 
  the generalized tree. 
  Note that $\overline{\Gamma}$ is not a simplex 
  because $\height I_{\Gamma} = \indeg I \geq 2$. 
  Then $\overline{J} := I_{\overline{\Gamma}}$ has a $2$-linear resolution.   

  \par
  If $\height \overline{J} = 1$, then $\overline{J}$ is of the form 
  $(x_1) \cap P_2$, and $I^{\ast}$ is of the form 
  \begin{displaymath}
    I^{\ast} = I_{\Gamma} 
             = P_{F \cup \{ x \} } \cap (x_1, x) \cap (P_2 R + (x)). 
  \end{displaymath}
  Therefore $\mu (I) \leq 3$. 

  \par
  Thus we may assume $\height \overline{J} \geq 2$. 
  Then $\overline{I} := (\overline{J})^{\ast}$ satisfies the assumptions 
  of Theorem \ref{claim:h2CM}. By the induction hypothesis, 
  we have $\ara \overline{I} = \pd_R R/{\overline{I}} = 2$. 
  Hence, we have $\ara I = \pd_R R/I = 2$ by Proposition \ref{claim:ara=2}. 
\end{proof}

\par
The next example, which is a generalization of Example \ref{exam:AdualLine}, 
gives an example of ideals which satisfy the assumptions of Theorem 
\ref{claim:h2CM}. 
\begin{example}
  \label{exam:AdualLineGen}
  Let us consider the squarefree monomial ideal $I_n$ 
  of $K[x_1, x_2, \ldots, x_n]$ 
  $(n \geq 4)$ generated by the following $n-1$ elements: 
  \begin{displaymath}
    m_i^{(n)} = \frac{x_1 \cdots x_n}{x_{n-i} x_{n-i+1}}, 
    \quad i = 1, 2, \ldots, n-1. 
  \end{displaymath}
  That is, $I_n$ is the Alexander dual ideal 
  of the Stanley--Reisner ideal $I_{{\Gamma}_n}$, 
  where ${\Gamma}_n$ is the simplicial complex whose facets are 
  $\{ x_1, x_2 \}, \{ x_2, x_3 \}, \ldots, \{ x_{n-1}, x_n \}$. 
  The ideals $I, I'$ in Example \ref{exam:AdualLine} are $I_4, I_5$, 
  respectively. 

  \par
  Then the height of $I_n$ is equal to $2$, and the quotient ring is 
  Cohen--Macaulay. Therefore by Theorem \ref{claim:h2CM}, we have 
  $\ara I_n = 2$. 

  \par
  For $n=4, 5$, two elements $q_1^{(n)}, q_2^{(n)}$ which generate $I_n$ 
  up to radical are given in Example \ref{exam:AdualLine}, i.e., 
  \begin{displaymath}
    \left\{ 
    \begin{aligned}
      q_1^{(4)} &= m_2^{(4)}, \\
      q_2^{(4)} &= m_1^{(4)} + m_3^{(4)}, 
    \end{aligned}
    \right. \qquad \left\{ 
    \begin{aligned}
      q_1^{(5)} &= x_5 q_1^{(4)} - x_1 x_2 m_1^{(5)}, \\
      q_2^{(5)} &= x_5 q_2^{(4)} - x_2 x_3 m_1^{(5)}. 
    \end{aligned}
    \right. 
  \end{displaymath}
  In general, two elements $q_1^{(n)}, q_2^{(n)}$ 
  which generate $I_n$ up to radical 
  are given by the following recursive formula: 
  \begin{equation}
    \label{eq:AdualLineGenerator}
    \left\{ 
    \begin{aligned}
      &q_1^{(n+1)} 
       = x_{n+1} q_1^{(n)} - x_{n-2}^{n-3} q_1^{(n-1)} m_1^{(n+1)}, \\
      &q_2^{(n+1)} 
       = x_{n+1} q_2^{(n)} - x_{n-2}^{n-3} q_2^{(n-1)} m_1^{(n+1)}, 
    \end{aligned}
    \right. 
    \qquad n \geq 5. 
  \end{equation}
%  We prove $\sqrt{(q_1^{(n)}, q_2^{(n)})} = I_n$ by induction on $n$. 
  We prove this by induction on $n$. 
  Note that $I_{n}' = I_{n+1}$ 
  with $F = \{ x_{n} \} (\subset G = \{ x_{n-1}, x_n \} )$ and $x_0 = x_{n+1}$ 
  with respect to the notations of the proof 
  of Proposition \ref{claim:ara=2}.  
  Hence by the proof of Proposition \ref{claim:ara=2}, 
  it suffices to check the following equality by induction on $n$ 
  under the hypothesis that $q_1^{(n)}, q_2^{(n)}$ generate 
  $I_{n}$ up to radical: 
%  under the hypothesis $\sqrt{(q_1^{(n)}, q_2^{(n)})} = I_{n}$: 
  \begin{equation}
    \label{eq:AdualLineGen}
    (m_1^{(n)})^{n-2} 
    = - x_{n-2}^{n-3} q_2^{(n-1)} q_1^{(n)} 
    + x_{n-2}^{n-3} q_1^{(n-1)} q_2^{(n)}, 
    \qquad n \geq 5. 
  \end{equation}
  When $n=5$, since
  \begin{displaymath}
    \begin{aligned}
        - q_2^{(4)} q_1^{(5)} + q_1^{(4)} q_2^{(5)} 
      &= - q_2^{(4)} (x_5 q_1^{(4)} - x_1 x_2 m_1^{(5)}) 
         + q_1^{(4)} (x_5 q_2^{(4)} - x_2 x_3 m_1^{(5)}) \\
      &= (x_1 x_2 q_2^{(4)} - x_2 x_3q_1^{(4)}) m_1^{(5)} \\
      &= (m_1^{(4)})^2 m_1^{(5)} \quad \text{by (\ref{eq:AdualLine})} 
    \end{aligned}
  \end{displaymath}
  and $x_3^2 (m_1^{(4)})^2 m_1^{(5)} = (m_1^{(5)})^3$, 
  we have the desired equality. 
  Similarly, for general $n$, 
  \begin{displaymath}
    \begin{aligned}
        &- q_2^{(n-1)} q_1^{(n)} 
         + q_1^{(n-1)} q_2^{(n)} \\
      = &- q_2^{(n-1)} 
         (x_n q_1^{(n-1)} - x_{n-3}^{n-4} q_1^{(n-2)} m_1^{(n)}) 
       + q_1^{(n-1)} 
         (x_n q_2^{(n-1)} - x_{n-3}^{n-4} q_2^{(n-2)} m_1^{(n)}) \\
      = &(x_{n-3}^{n-4} q_1^{(n-2)} q_2^{(n-1)} 
          - x_{n-3}^{n-4} q_2^{(n-2)} q_1^{(n-1)}) m_1^{(n)} \\
      = &(m_1^{(n-1)})^{n-3} m_1^{(n)} \quad \text{by induction hypothesis} 
    \end{aligned}
  \end{displaymath}
  and $x_{n-2}^{n-3} (m_1^{(n-1)})^{n-3} m_1^{(n)} = (m_1^{(n)})^{n-2}$ yield 
  the equation (\ref{eq:AdualLineGen}). 
\end{example}

\par
Another class of ideals which satisfies the assumptions of Theorem 
\ref{claim:h2CM} is found in Barile \cite[Section 3]{Barile96}. 
It is essentially Alexander dual of the class of 
Ferrers ideals (see \cite{Corso-Nagel}, \cite{Barile0606}). 
In \cite{Barile96}, Barile construct $2$ elements which generate the ideals 
up to radical 
%proved $\ara I = \pd_R R/I = 2$ 
on the different way.

\section{Improvement proof of the result by Barile and Terai}
Let $\Delta$ be a simplicial complex on the vertex set $X$. 
Let $F$ be a face of $\Delta$ and $x_0$ new vertex. 
Set ${\Delta}' = \Delta \cup \cone_{x_0} F$. 
Throughout of this section, we use these notations. 
Note that these are different from those of previous sections. 

\par
In our motivated paper Barile and Terai \cite{BariTera08}, 
the main result \cite[Theorem 1]{BariTera08} 
depends on the base field $K$. 
Precisely, it needs the assumption that $K$ is algebraically closed. 
In this section, we give an improved proof of it 
which does not depend on the base field $K$. 
\begin{theorem}[{cf.\  Barile and Terai \cite[Theorem 1]{BariTera08}}]
  \label{claim:Mthm}
  Let $\Delta$ be a simplicial complex on the vertex set 
  $X = \{ x_1, x_2, \ldots, x_n \}$, 
  $F$ a face of $\Delta$, and $x_0$ new vertex. 
  Set ${\Delta}' = \Delta \cup \cone_{x_0} F$. Then 
  \begin{displaymath}
    \ara I_{{\Delta}'} \leq \max \{ \ara I_{\Delta} + 1, n - |F| \}. 
  \end{displaymath}
%  where $|F|$ denotes the cardinality of $F$. 
\end{theorem}

\par
As a consequence of our improvement, 
we can also omit the assumption on $K$ 
for other results in \cite{BariTera08}: 
\begin{theorem}[{cf.\  Barile and Terai \cite[Theorem 2]{BariTera08}}]
  \label{claim:cone-thm}
  Let $\Delta$ be a simplicial complex on the vertex set 
  $X = \{ x_1, x_2, \ldots, x_n \}$, 
  $F$ a face of $\Delta$, and $x_0$ new vertex. 
  Set ${\Delta}' = \Delta \cup \cone_{x_0} F$. 
  If $\ara I_{\Delta} = \pd K[{\Delta}]$ holds, then 
  $\ara I_{{\Delta}'} = \pd K[{\Delta}']$ also holds. 
\end{theorem}

\begin{corollary}[{cf.\  Barile and Terai \cite[Corollary 3]{BariTera08}}]
  \label{claim:2-linear}
  Let $I$ be a squarefree monomial ideal of $R=K[X]$. 
  Suppose that $I$ has a $2$-linear resolution. 
  Then 
  \begin{displaymath}
    \ara I = \pd_R R/I. 
  \end{displaymath}
\end{corollary}
Corollary \ref{claim:2-linear} was first proved 
by Morales \cite[Theorems 8 and 9]{Mo} 
on the different way, but he also assumed that $K$ is algebraically closed. 

\par
\bigskip

\par
Now, we prove Theorem \ref{claim:Mthm}. 
The proof is divided into two steps. 
We construct $\max \{ \ara I_{\Delta} + 1, n-|F| \}$ elements 
which generate $I_{{\Delta}'}$ up to radical in the latter step (Step 2). 
The former step (Step 1) is assigned to transform elements 
which generate $I_{\Delta}$ up to radical 
so that the elements constructed in (Step 2) belong to $I_{{\Delta}'}$. 

\par
In our proof, (Step 1) is the same as that by Barile and Terai
(see also Barile \cite[Theorem 1]{Barile08-1}). 
Thus we omit the detail. 
Our improvement is in (Step 2). 
In Case 1 of (Step 2), 
the elements which generate $I_{{\Delta}'}$ up to radical 
are the same as those of Barile and Terai. 
The difference is that we use the cofactor matrix 
instead of Cramer's Rule which they used, 
and we do not use Hilbert's Nullstellensatz. 
%(The way of proof using cofactor matrix is also useful 
%for \cite[Theorem 1]{Barile08-1}.) 
In Case 2 of (Step 2), we give the different elements 
which generate $I_{{\Delta}'}$ up to radical from those of Barile and Terai. 
This is our main improvement. 

\begin{proof}[Proof of Theorem \ref{claim:Mthm}]
  (\textit{Step $1$}) 
  First, we fix the notation. 
  Set $R = K[X]$ and $R' = K[X']$ where $X' = X \cup \{ x_0 \}$. 
  If $F=X$, then $I_{\Delta} = I_{{\Delta}'}=0$ and the assertion is trivially 
  true. Thus we assume $F \neq X$. 
  Let $G$ be a facet of $\Delta$ which contains $F$. 
  We can assume that $G = \{ x_{s+1}, \ldots, x_n \}$ 
  and $F = \{ x_{t+1}, \ldots, x_n \}$, where $s \leq t$. 
  Then $I_{{\Delta}'} = I_{\Delta} R' + (x_0 x_1, \ldots, x_0 x_t)R'$. 
  We set $\ara I_{\Delta} = h$. 
  Then we can rewrite the claim as 
  \begin{displaymath}
    \ara I_{{\Delta}'} \leq \max \{ h+1, t \}. 
  \end{displaymath}

  \par
  Assume that $q_1, \ldots, q_h$ generate $I_{\Delta}$ up to radical. 
  Since $I_{\Delta} \subset P_{G} = (x_1, \ldots, x_s)$ 
  and $q_i \in I_{\Delta}$, 
  we can write 
  \begin{displaymath}
    q_i = \sum_{j=1}^s a_{ij} x_{j}, \qquad i = 1, 2, \ldots, h, 
  \end{displaymath}
  where $a_{ij} \in R$. 
  Then we can transform each $q_i$ to 
  \begin{displaymath}
    \overline{q}_i = \sum_{j=1}^s \overline{a}_{ij} x_j, 
  \end{displaymath}
  where $\overline{a}_{ij} \in I_{\Delta}$ preserving the property that 
  $\overline{q}_1, \ldots, \overline{q}_h$ generate $I_{\Delta}$ 
  up to radical; see Barile and Terai \cite[Proof of Theorem 1]{BariTera08}. 

  \par
  (\textit{Step $2$}) 
  Now we find $\max \{ h+1, t \}$ elements which generate $I_{{\Delta}'}$ 
  up to radical. We divide into two cases. 

  \begin{description}
    \item[Case 1] Suppose that $h+1 > t$. 
  \end{description}
  We show that $\ara I_{{\Delta}'} \leq h+1$. 
  We set $\overline{A} = (\overline{a}_{ij})_{i,j = 1, \ldots, t}$, 
  where $\overline{a}_{ij} = 0$ if $j>s$. 
  Let $A_1 = \overline{A} + x_0 Id_t$, where $Id_t$ denotes the $t \times t$ 
  identity matrix. Set 
  \begin{displaymath}
    J_1 = (\det A_1 - x_0^t, 
           \overline{q}_1 + x_0 x_1, \ldots, \overline{q}_t + x_0 x_t, 
           \overline{q}_{t+1}, \ldots, \overline{q}_h) R'. 
  \end{displaymath}
  We prove $\sqrt{J_1} = I_{{\Delta}'}$. 
  Since $\overline{a}_{ij} \in I_{\Delta}$, 
  we have $\det A_1 - x_0^t \in I_{\Delta} R'$. 
  Moreover since $\overline{q}_i \in I_{\Delta}$, $i= 1, 2, \ldots, h$ and 
  $x_0 x_j \in I_{{\Delta}'}$, $j=1, 2, \ldots, t$, 
  we have $\sqrt{J_1} \subset I_{{\Delta}'}$. 
  We prove the opposite inclusion. 
  To do this, it suffices to show that $\overline{q}_i \in \sqrt{J_1}$, 
  $i=1, 2, \ldots, t$ and $x_0 x_j \in \sqrt{J_1}$, $j=1, 2, \ldots, t$. 

  \par
  Let $B_1$ be the cofactor matrix of $A_1$. Then $B_1 A_1 = (\det A_1) Id_t$. 
  Since 
  \begin{displaymath}
    \left(
    \begin{matrix}
      \overline{q}_1 + x_0 x_1 \\
      \vdots \\
      \overline{q}_t + x_0 x_t 
    \end{matrix}
    \right) = A_1 \left(
    \begin{matrix}
      x_1 \\
      \vdots \\
      x_t
    \end{matrix}
    \right), 
  \end{displaymath}
  we have 
  \begin{displaymath}
    B_1 \left(
    \begin{matrix}
      \overline{q}_1 + x_0 x_1 \\
      \vdots \\
      \overline{q}_t + x_0 x_t 
    \end{matrix}
    \right) = B_1 A_1 \left(
    \begin{matrix}
      x_1 \\
      \vdots \\
      x_t
    \end{matrix}
    \right) = (\det A_1) \left( 
    \begin{matrix}
      x_1 \\
      \vdots \\
      x_t
    \end{matrix}
    \right). 
  \end{displaymath}
  Then $(\det A_1) x_j \in J_1$ for $j = 1, 2, \ldots, t$ 
  by $\overline{q}_i + x_0 x_i \in J_1$ for all $i=1, 2, \ldots, t$. 
  Multiplying $\det A_1 - x_0^t \in J_1$ by $x_j$, 
  we have $x_0^t x_j \in J_1$. 
  Hence $x_0 x_j \in \sqrt{J_1}$ for $j=1, 2, \ldots, t$. 
  Since $\overline{q}_i + x_0 x_i \in J_1$, 
  we have $\overline{q}_i \in \sqrt{J_1}$ 
  for $i=1, 2, \ldots, t$, as required. 

  \begin{description}
    \item[Case 2] Suppose that $h+1 \leq t$. 
  \end{description}
  We show that $\ara I_{{\Delta}'} \leq t$. 
  Note that in this case, $s \leq t-1$ because if $s=t$, then 
  $t$ is the height of the minimal prime $P_{G}$ of $I_{\Delta}$ and 
  Krull's principal ideal theorem shows that $\ara I \geq t$. 
%  \begin{displaymath}
%    \ara I_{\Delta} \geq \pd K[{\Delta}] \geq \bight I_{\Delta} \geq t, 
%  \end{displaymath}
%  where $\bight I_{\Delta}$ denotes maximum height of minimal prime divisors 
%  of $I_{\Delta}$. 
  This contradicts to $\ara I_{\Delta} = h \leq t-1$. 
%  Therefore $\overline{a}_{it} = 0$ for all $i = 1, 2, \ldots, h$. 

  \par
  We set $\overline{A}' = (\overline{a}_{ij})_{i,j = 1, \ldots, t-1}$, 
  where $\overline{a}_{ij} = 0$ if $i>h$ or $j>s$. 
  Let $A_2 = \overline{A}' + x_0 Id_{t-1}$, where $Id_{t-1}$ denotes 
  the $(t-1) \times (t-1)$ identity matrix. Set 
  \begin{displaymath}
    J_2 = ((\det A_2)(x_0+x_t) - x_0^t, 
           \overline{q}_1 + x_0 x_1, \ldots, \overline{q}_h + x_0 x_h, 
           x_0 x_{h+1}, \ldots, x_0 x_{t-1}) R'. 
  \end{displaymath}
  We prove $\sqrt{J_2} = I_{{\Delta}'}$. 
  As $\overline{a}_{ij} \in I_{\Delta}$, similarly to Case 1, 
  we have $\sqrt{J_2} \subset I_{{\Delta}'}$. 
  We prove the opposite inclusion. 
  Let $B_2$ be the cofactor matrix of $A_2$. 
  Then $B_2 A_2 = (\det A_2) Id_{t-1}$. 
  Since we set $\overline{a}_{ij} = 0$ for $i>h$, we can write formally 
  $x_0 x_i = \overline{q}_i + x_0 x_i$ for $i = h+1, \ldots, t-1$. 
  Using this notation, we have 
  \begin{displaymath}
    \left(
    \begin{matrix}
      \overline{q}_1 + x_0 x_1 \\
      \vdots \\
      \overline{q}_{t-1} + x_0 x_{t-1} 
    \end{matrix}
    \right) = A_2 \left(
    \begin{matrix}
      x_1 \\
      \vdots \\
      x_{t-1}
    \end{matrix}
    \right). 
  \end{displaymath}
  Thus 
  \begin{displaymath}
    B_2 \left(
    \begin{matrix}
      \overline{q}_1 + x_0 x_1 \\
      \vdots \\
      \overline{q}_{t-1} + x_0 x_{t-1} 
    \end{matrix}
    \right) = B_2 A_2 \left(
    \begin{matrix}
      x_1 \\
      \vdots \\
      x_{t-1}
    \end{matrix}
    \right) = (\det A_2) \left( 
    \begin{matrix}
      x_1 \\
      \vdots \\
      x_{t-1}
    \end{matrix}
    \right). 
  \end{displaymath}
  Then $(\det A_2) x_j \in J_2$ for $j = 1, 2, \ldots, t-1$ 
  by $\overline{q}_i + x_0 x_i \in J_2$ for all $i=1, 2, \ldots, t-1$. 
  Multiplying $(\det A_2)(x_0+x_t) - x_0^t \in J_2$ by $x_j$, 
  we have $x_0^t x_j \in J_2$ for $j = 1, 2, \ldots, t-1$. 
  Hence $x_0 x_j \in \sqrt{J_2}$ for $j=1, 2, \ldots, t-1$. 
  Since $\overline{q}_i + x_0 x_i \in J_2$, 
  we have $\overline{q}_i \in \sqrt{J_2}$ for $i=1, 2, \ldots, t-1$. 
  In particular, $\overline{q}_i \in \sqrt{J_2}$ for $i=1, 2, \ldots, h$. 
  Since $\sqrt{(\overline{q}_1, \ldots, \overline{q}_h)} = I_{\Delta}$ 
  and $\overline{a}_{ij} \in I_{\Delta}$, 
  we have 
  \begin{displaymath}
    \overline{a}_{ij} \in \sqrt{(\overline{q}_1, \ldots, \overline{q}_h)} 
    \subset \sqrt{J_2}, 
    \qquad \text{for all $i, j$}. 
  \end{displaymath}
  Therefore $(\det A_2)(x_0 + x_t) - x_0^t \in J_2$ implies 
  $x_0^{t-1} (x_0 + x_t) - x_0^t \in \sqrt{J_2}$. 
  Thus we have $x_0^{t-1} x_t \in \sqrt{J_2}$, 
  that is $x_0 x_t \in \sqrt{J_2}$. 
  This completes the proof. 
\end{proof}

The next example was considered in Barile and Terai 
\cite[Example 1]{BariTera08} as the example to show elements 
given in the proof of \cite[Theorem 1]{BariTera08}. 
We show the elements given in our proof at the same example, too. 
\begin{example}
  \label{ex:cone}
  Let $\Delta$ be the simplicial complex on the vertex set 
  $\{ x_1, x_2, x_3, x_4 \}$ whose facets are 
  $\{ x_1, x_2 \}, \{ x_1, x_4 \}, \{ x_2, x_3 \}, \{ x_3, x_4 \}$. 
  Then 
  \begin{displaymath}
    I_{\Delta} = (x_1, x_2) \cap (x_1, x_4) \cap (x_2, x_3) \cap (x_3, x_4) 
               = (x_1 x_3, x_2 x_4). 
  \end{displaymath}
  Thus $I_{\Delta}$ is complete intersection. In particular, 
  $h = \ara I_{\Delta} = 2$. Set $q_1 = x_1 x_3$ and $q_2 = x_2 x_4$. 
  Take the face $F = \{ x_4 \} \in \Delta$ and let $x_0$ be a new vertex. 
  Then $I_{{\Delta}'}$ is generated by the following $5$ elements: 
  \begin{displaymath}
    x_1 x_3, x_2 x_4, x_0 x_1, x_0 x_2, x_0 x_3. 
  \end{displaymath}
  Then $t=3$. 
  We take the facet $G$ as $\{ x_3, x_4 \}$. 
  Then $P_{G} = (x_1, x_2)$. 
  In this case, we have 
  \begin{displaymath}
    \overline{q}_1 = x_1 x_3^2 \cdot x_1, 
    \qquad \overline{q}_2 = x_2 x_4^2 \cdot x_2. 
  \end{displaymath}
  Since $h+1=3=t$, 
  we apply the Case 2 of the proof of Theorem \ref{claim:Mthm}. 
  Since 
  \begin{displaymath}
    \overline{A}' = \left(
    \begin{matrix}
      x_1 x_3^2 & 0 \\
      0 & x_2 x_4^2
    \end{matrix}
    \right), \qquad A_2 = \left(
    \begin{matrix}
      x_1 x_3^2 + x_0 & 0 \\
      0 & x_2 x_4^2 + x_0
    \end{matrix}
    \right), 
  \end{displaymath}
  the ideal $I_{{\Delta}'}$ is generated by the following $3$ elements 
  up to radical: 
  \begin{displaymath}
    \begin{aligned}
      &(x_1 x_3^2 + x_0)(x_2 x_4^2 + x_0)(x_0 + x_3) - x_0^3 \\
      &= x_0^2 x_3 
       + x_0^2 x_1 x_3^2+ x_0^2 x_2 x_4^2 + x_0 x_1 x_2 x_3^2 x_4^2 
       +x_0 x_1 x_3^3 + x_0 x_2 x_3 x_4^2 + x_1 x_2 x_3^3 x_4^2, \\
      &x_1^2 x_3^2 + x_0 x_1, \\
      &x_2^2 x_4^2 + x_0 x_2. 
    \end{aligned}
  \end{displaymath}
\end{example}

%%%%%%%%%%%%%%%%%%%%%%%%%%%%%%%%%%%%%%%%
%%% Acknowledgments
\begin{acknowledgements}
  The author is grateful to Naoki Terai for drawing the author's attention 
  to this problem. 
  The author also thanks to Ken-ichi Yoshida for giving her 
  valuable suggestions. 

  \par
  This research was supported by JST, CREST. 
\end{acknowledgements}

%%%%%%%%%%%%%%%%%%%%%%%%%%%%%%%%%%%%%%%%
%%%%%%%%%%%%%%%%%%%%%%%%%%%%%%%%%%%%%%%%
%%%%%%%%%%%%%%%%%%%%%%%%%%%%%%%%%%%%%%%%
%%% references
%\bibliographystyle{amsplain}

\end{document}